\documentclass[11pt]{article}

\usepackage{amsmath}    
\usepackage{graphicx}   
\usepackage{verbatim}   
\usepackage{color}      
\usepackage{subfigure}  
\usepackage{amssymb,tikz}
\usepackage{soul}

\newtheorem{theorem}{Theorem}[section]
\newtheorem{proposition}[theorem]{Proposition}

\newtheorem{lemma}[theorem]{Lemma}
\newtheorem{remark}[theorem]{Remark}

\newcommand{\proof}{\noindent{\bf Proof.\ }}
\newcommand{\qed}{\hfill $\square$\medskip}

\newcommand{\R}{\mathbb R}
\newcommand{\N}{\mathbb N}
\newcommand{\Z}{\mathbb Z}

\newcommand{\V}{\mathcal V}

\DeclareMathOperator {\g} {g}
\DeclareMathOperator {\sg} {sg}
\DeclareMathOperator {\ip} {ip}
\DeclareMathOperator {\diam} {diam}
\DeclareMathOperator {\rad} {rad}
\let\d\relax
\DeclareMathOperator {\d} {d}
\DeclareMathOperator {\n} {n}

\begin{document}

\tikzstyle{every node}=[circle, draw, fill=black!10,
                        inner sep=0pt, minimum width=4pt]

\title{Strong geodetic number of complete bipartite graphs and of graphs with specified diameter }
 
\author{
    Vesna Ir\v si\v c $^a$
}

\date{July 20, 2017}

\maketitle
\begin{center}

	$^a$ Institute of Mathematics, Physics and Mechanics, Ljubljana, Slovenia\\
	{\tt vesna.irsic@student.fmf.uni-lj.si}\\
	\medskip

\end{center}

\begin{abstract}
The strong geodetic problem is a recent variation of the classical geodetic problem. For a graph $G$, its strong geodetic number $\sg(G)$ is the cardinality of a smallest vertex subset $S$, such that each vertex of $G$ lies on one fixed geodesic between a pair of vertices from $S$. In this paper, some general properties of the strong geodesic problem are studied, especially in connection with diameter of a graph. The problem is also solved for balanced complete bipartite graphs. 
\end{abstract}

\noindent{\bf Keywords:} geodetic number; strong geodetic number; isometric path number;   complete bipartite graphs; diameter

\medskip
\noindent{\bf AMS Subj.\ Class.: 05C12, 05C70}

\section{Introduction}

Covering vertices of a graph with shortest paths is a problem with several different variations. The {\em geodetic number}~\cite{HLT93} of a graph $G$, denoted by $\g(G)$, is the cardinality of the smallest subset of vertices such that the geodesics between them cover all vertices. For results up to 2011 see~\cite{BrKo11}. Applications of geodetic problem can be found in convexity theory~\cite{Centeno, Jiang+2004, Lu2004, pelayo-2013} and game theory~\cite{igre}. Connection between geodetic number and diameter of a graph has been studied in~\cite{obstoj}. See~\cite{Filomat} for characterization of graphs with large geodetic number and~\cite{Scientia0, Scientia} for results on the edge geodetic problem. For additional results see~\cite{ekim-2014, ekim-2012} and~\cite{soloff-2015}.

Another variation is the {\em isometric path number}~\cite{fisher-2001} which is defined as the minimum number of shortest paths required to cover all vertices of a graph $G$ and is denoted by $\ip(G)$. See~\cite{Fitzpatrick-1999-CN} for results about isometric path number of the Cartesian products of graphs and~\cite{pan-2006} for exact values for complete $r$-partite graphs and for some Hamming graphs.

The strong geodetic problem was introduced in~\cite{MaKl16a} as another variation of the problem. It was motivated by applications in social networks and is defined as follows. Let $G=(V,E)$ be a graph. Given a set $S\subseteq V$, for each pair of vertices $\{x,y\}\subseteq S$, $x\ne y$, let $\widetilde{g}(x,y)$ be a {\em selected fixed} shortest path between $x$ and $y$. We set  
$$\widetilde{I}(S)=\{\widetilde{g}(x, y) : x, y\in S\}\,,$$ and $V(\widetilde{I}(S))=\bigcup_{\widetilde{P} \in \widetilde{I}(S)} V(\widetilde{P})$. If $V(\widetilde{I}(S)) = V$ for some $\widetilde{I}(S)$, then the set $S$ is called a {\em strong geodetic set}. For a graph $G$ with just one vertex, we take $V(G)$ as its unique strong geodetic set. The {\em strong geodetic problem} is to find a minimum strong geodetic set $S$ of $G$. Clearly, the collection $\widetilde{I}(S)$ of geodesics consists of exactly $\binom{|S|}{2}$ paths. The cardinality of a minimum strong geodetic set is the {\em strong geodetic number} of $G$ and is denoted by $\sg(G)$. See also~\cite{Klavzar+2017} for additional results on the strong geodetic number and~\cite{MaKl16b} for an edge version of the problem.

The rest of the paper is organized as follows. In Section~\ref{sec:bipartite} we determine the strong geodetic number of balanced complete bipartite graphs while in Section~\ref{sec:diameter} we consider the connection between strong geodetic number and diameter. We conclude the paper with some suggestions for further investigation. But first we define concepts needed. 

All graphs considered in this paper are simple and connected. The {\em distance} $d_G(u,v)$ between vertices $u$ and $v$ of a graph $G$ is the number of edges on a shortest $u,v$-path alias $u,v$-{\em geodesic}. The {\em diameter} ${\rm diam}(G)$ of $G$ is the maximum distance between the vertices of $G$. The {\em radius} $\rad(G)$ of $G$ is the minimum eccentricity of all vertices, i.e.\ $\rad(G) = \min_{u \in V(G)} \max_{v \in V(G)} \d(u, v)$. We denote the order of a graph by $\n(G)$. A vertex $v$ of a graph $G$ is {\em simplicial} if its neighborhood induces a clique. We will use the notation $[n] = \{1,\ldots, n\}$ and the convention that $V(P_n) = [n]$ for any $n\ge 1$, where the edges of $P_n$ are defined in the natural way.

To conclude the introduction we state the following simple but fundamental result and recall the observation from~\cite{MaKl16a} which asserts that a simplicial vertex necessarily lies in any strong geodetic set of a graph.

\begin{proposition}
\label{prop:basic}
If $G$ is a graph with $\n(G) = n \geq 2$, then $$2 \leq \sg(G) \leq n.$$
Moreover, $\sg(G) = 2$ if and only if $G \cong P_n$, and $\sg(G) = n$ if and only if $G \cong K_n$.
\end{proposition}

\proof
The inequality $2 \leq \sg(G) \leq n$ is clear (having in mind that $n \geq 2$).  

If $G \cong P_n$, then $\sg(G) = 2$ as $\{1, n\}$ is a strong geodetic set. Now let $G$ be a graph with $n$ vertices and $\sg(G) = 2$. Then there exist two vertices, say $u, v \in V(G)$, such that $S = \{u, v\}$ is a strong geodetic set. All other vertices of the graph $G$ must therefore lie on a fixed $u, v$-geodesic. Hence graph $G$ consists only of this geodesic (if there were other edges, the distance between $u$ and $v$ would change). Thus we have $G \cong P_n$. 

If $G \cong K_n$, then all vertices of $G$ are simplicial and must therefore lie in any strong geodetic set. Hence $\sg(G) = n$. Now let $G$ be a graph with $n$ vertices and $\sg(G) = n$. Suppose there exist $u, v \in V(G)$ which are not neighbours. Then $u, v$-geodesic has length at least $2$, so there exists some $w \in V(G) - \{u, v\}$ on it. Consider the set $S = V(G) - \{w\}$. Geodesics between $u$ and other vertices in $S -  \{v\}$ cover all vertices in $V(G) - \{v, w\}$ and the $u, v$-geodesic covers vertices $v$ and $w$. Hence $S$ is a strong geodetic set and $\sg(G) \leq n-1$ which contradicts $\sg(G) = n$.
\qed

\section{Complete bipartite graphs}
\label{sec:bipartite}

In this section we study the strong geodetic problem on complete bipartite graphs and determine the exact values for graphs $K_{n,n}$ and $K_{n_1, n_2}$ where $n_1 \gg n_2$. In the case when $\min\{m,n\} = 1$, the strong geodetic problem is trivial. Therefore we only consider the case $m, n \geq 2$. 

The geodetic number of complete bipartite graphs~\cite{HLT93} is $$\g(K_{m, n}) = \min\{ m, n, 4 \}$$  and the isometric path number~\cite{pan-2006} equals $$\ip(K_{m, n}) = \begin{cases}
\left \lceil \frac{\max\{m, n\}}{2} \right \rceil; & \max\{m, n\} > 2 \min\{m, n\},\\
\left \lceil \frac{m + n}{3} \right \rceil; & \text{otherwise}.
\end{cases}$$ 

On the other hand, the strong geodetic problem on complete bipartite graphs is more complex and has not been studied before. 

Every geodesic in a complete bipartite graph is either an edge or a path of length $2$ with both endvertices in the same part of the bipartition. If a strong geodetic set $S$ has $k$ vertices in one part of the bipartition, then geodesics between those vertices can cover at most $\binom{k}{2}$ vertices in the other part (as for each pair of vertices from $S$ we can fix a geodesic through one vertex in the other part). Hence the strong geodetic problem for complete bipartite graph can be transformed into a (nonlinear) integer program. But first we need to specify some notation.

Let $(X, Y)$ be the bipartition of $K_{n_1, n_2}$ and $S = S_1 \cup S_2$, $S_1 \subseteq X$, $S_2 \subseteq Y$, its strong geodetic set.  Let $|S_i| = s_i$ for $i \in [2]$. Thus $\sg(K_{n_1, n_2}) = s_1 + s_2$. With geodesics between vertices from $S_1$ we wish to cover vertices in $Y - S_2$. And vice versa, with geodesics between vertices from $S_2$ we are covering vertices in $X - S_1$. The optimization problem for $\sg(K_{n_1, n_2})$ now reads as follows:
\begin{align}
\label{optimizationProblem}
\begin{split}
\min \quad & s_1 + s_2 \\
\text{subject to: } & 0 \leq s_1 \leq n_1\\
&  0 \leq s_2 \leq n_2\\
&  \binom{s_2}{2} \geq n_1 - s_1\\
&  \binom{s_1}{2} \geq n_2 - s_2\\
& s_1, s_2 \in \Z.
\end{split}
\end{align}

In the rest of the section we consider two special cases, $n_1 = n_2$ and $n_1 \gg n_2$,  and determine the strong geodetic number for them.

First we focus on the case when $n_1 = n_2 = n$. By solving concrete optimization problems, we get $\sg(K_{2,2}) = \sg(K_{3,3}) = 3$, $\sg(K_{4,4}) = 5$ and $\sg(K_{5,5}) = 5$. For larger values of $n$ we have the following formula which was conjectured, using a computer experiment, by Petkov\v sek~\cite{petkovsek}.

\begin{theorem}
\label{sg polni dvodelni}
If $n \geq 6$, then
$$\sg(K_{n,n}) = \begin{cases}
2 \left \lceil \displaystyle \frac{-1 + \sqrt{8 n + 1}}{2} \right \rceil; & 8n - 7 \text{ is not a perfect square},\vspace{0.2cm}\\
2 \left \lceil \displaystyle \frac{-1 + \sqrt{8 n + 1}}{2} \right \rceil - 1; & 8n - 7 \text{ is a perfect square}.
\end{cases}$$
For an optimal strong geodetic set $S = S_1 \cup S_2$ it holds $s_1 = s_2 = \left \lceil \frac{-1 + \sqrt{8 n + 1}}{2} \right \rceil$, if $8n - 7$ is not a perfect square, otherwise we have $s_1 = \left \lceil \frac{-1 + \sqrt{8 n + 1}}{2} \right \rceil$ and $s_2 = s_1 - 1$.
\end{theorem}

Intuitively, the optimal value is reached when the strong geodetic set contains approximately  the same number of vertices in each part of the bipartition. This is formally shown below and motivates the next technical lemma.

\begin{lemma}
\label{lema sg polni dvodelni}
Let $T = T_1 \cup T_2$ be a strong geodetic set of $K_{n,n}$, $n \geq 6$, with bipartition $(X, Y)$, where $T_1 \subseteq X$, $T_2 \subseteq Y$ and $t_i = |T_i|$ for all $i \in [2]$. If $|t_1 - t_2| \geq 2$, then there exists a strong geodetic set $T' = T'_1 \cup T'_2$, $T'_1 \subseteq X$, $T'_2 \subseteq Y$, such that $|T'| = |T|$ and $|t'_1 - t'_2| < |t_1 - t_2|$, where $t'_i = |T'_i|$ for $i \in [2]$.
\end{lemma}

\proof
Without loss of generality we assume $t_1 \geq t_2$. Let $t_1 - t_2 = k \geq 2$.

Assume first that $\min\{t_1, t_2\} \geq 1$. Let $V(K_{n,n}) = \{x_1, \ldots, x_n\} \cup \{y_1, \ldots, y_n\}$, where $x_i \in X$ and $y_i \in Y$ for all $i \in [n]$. Without loss of generality we assume $T_1 = \{x_1, \ldots, x_{t_2 + k}\}$ and $T_2 = \{ y_1, \ldots, y_{t_2} \}$. 

As $T$ is a strong geodetic set,
\begin{equation}
\label{eq1}
t_2 + \binom{t_2 + k}{2} \geq n
\end{equation}
and
\begin{equation}
\label{eq2}
t_2 + k + \binom{t_2}{2} \geq n.
\end{equation}

Define the set $T' = T'_1 \cup T'_2$, where $T'_1 = T_1 - \{x_{t_2 + k}\}$ and $T'_2 = T_2 \cup \{ y_{t_2 + 1} \}$. From the assumptions of the lemma it follows that $0 \leq t_i \leq n$ for $i \in [2]$. We now prove that $T'$ is a strong geodetic set. 

As $t_2 \geq 1$ we have $\binom{t_2 + 1}{2} \geq \binom{t_2}{2} + 1$. From this and~\eqref{eq2} it follows that $$t_2 + k - 1 + \binom{t_2 + 1}{2} \geq n.$$ 

It now suffices to show that all the vertices in $Y - T'_2$ can be covered with geodesics between vertices in $T'_1$. From~\eqref{eq2} it follows that with geodesics between $t_2$ vertices $\{x_1, \ldots, x_{t_2}\}$ we can cover $n - t_2 - k$ vertices $\{ y_{t_2 + k + 1}, \ldots, y_n \}$. 

We only need to prove that using the remaining $k-1$ vertices from the set $T'_1$ we can cover the remaining $k-1$ vertices in $Y - T'_2$. This clearly holds if $k-1 \geq 3$ i.e.\ $k \geq 4$, as in this case $\binom{k-1}{2} \geq k-1$. We can also cover the remaining vertices with geodesics between vertex $x_1$ and vertices from $\{ x_{t_2 + 1}, \ldots, x_{t_2 + k - 1} \}$ whenever $t_2 \geq k-1$. The only unsolved case is $k = 3$ and $t_2 = 1$, which is impossible as $n > 4$.  

Hence $T'$ is a strong geodetic set, such that $$|t'_1 - t'_2| = (t_1 - 1) - (t_2 + 1) = |t_1 - t_2| - 2 < |t_1 - t_2|.$$ 

Now consider the case $\min\{t_1, t_2\} < 1$. This means $t_2 = 0$. In this case $T$ can be a strong geodetic set if and only if $t_1 = n$. Define the set $T' = T'_1 \cup T'_2$, where $T'_1 = \{ x_1, \ldots, x_{\left \lceil \frac{n}{2} \right \rceil} \}$ and $T'_2 = \{ y_1, \ldots, y_{\left \lfloor \frac{n}{2} \right \rfloor} \}$. It holds $T'_1 \subseteq X$ and $T'_2 \subseteq Y$. From $n \geq 6$ it follows $\left \lceil \frac{n}{2} \right \rceil, \left \lfloor \frac{n}{2} \right \rfloor \geq 3$, so for those two integers we have $\binom{k}{2} \geq k$. Thus $T'$ is a strong geodetic set, such that $|t'_1 - t'_2| \in \{0, 1\}$ hence, $|t'_1 - t'_2| < |t_1 - t_2| = n$.
\qed

For the proof of Theorem~\ref{sg polni dvodelni} we need the following fact.

\begin{lemma}
\label{lihi}
For every odd perfect square $s$ there exists an integer $k$ such that $s = 8 k + 1$.
\end{lemma}

\proof
As $s$ is an odd perfect square, there exists an integer $l$ such that $s = (2 l + 1)^2 = 4 l (l + 1) + 1$. As $l (l+1)$ is even, we have $s = 8 l' + 1$. 
\qed

Now we can prove the main result of this section.

\proof[of Theorem~\ref{sg polni dvodelni}]
It follows from Lemma~\ref{lema sg polni dvodelni} that we only need to study strong geodetic sets $S = S_1 \cup S_2$ for which $|s_1 - s_2| \leq 1$. This gives rise to two cases. 

\begin{enumerate}
\item Suppose $s_2 = s_1 - 1$. The optimization problem~\eqref{optimizationProblem}  simplifies to
\begin{align*}
\min \quad & 2 s_1 - 1 \\
\text{subject to: } & 0 \leq s_1 \leq n\\
& s_1^2 - s_1 + 2 - 2n \geq 0\\
& s_1^2 + s_1 - 2 - 2n \geq 0\\
& s_1 \in \Z.
\end{align*}

It follows from the first quadratic inequality that $s_1 \geq \frac{1 + \sqrt{8 n - 7}}{2}$ and from the second that $s_1 \geq \frac{- 1 + \sqrt{8n + 9}}{2}$. As for all integers $n \geq 2$ it holds 
$$n \geq \frac{1 + \sqrt{8 n - 7}}{2} \geq \frac{- 1 + \sqrt{8n + 9}}{2},$$
minimum of $2 s_1 - 1$ is reached at $$s_1 = \left \lceil \frac{1 + \sqrt{8 n - 7}}{2} \right \rceil$$ and equals $2 \left \lceil \frac{1 + \sqrt{8 n - 7}}{2} \right \rceil - 1$.

\item Suppose $s_1 = s_2$. Then~\eqref{optimizationProblem} simplifies to
\begin{align*}
\min \quad & 2 s_1 \\
\text{subject to: } & 0 \leq s_1 \leq n_1\\
& s_1^2 + s_1 - 2n \geq 0\\
& s_1 \in \Z.
\end{align*}

Quadratic inequality yields $s_1 \geq \left \lceil \frac{-1 + \sqrt{8 n + 1}}{2} \right \rceil$, thus the minimum of $2 s_1$ is reached at
$$s_1 = \left \lceil \frac{-1 + \sqrt{8 n + 1}}{2} \right \rceil$$ and equals $2 \left \lceil \frac{-1 + \sqrt{8 n + 1}}{2} \right \rceil$.

\end{enumerate}

We must now consider which case gives rise to a smaller value of $s_1 + s_2$ and is thus equal to the strong geodetic number. Let $a = \left \lceil \frac{1 + \sqrt{8 n - 7}}{2} \right \rceil$, $b = \left \lceil \frac{-1 + \sqrt{8 n + 1}}{2} \right \rceil$, $\sg_a = 2 a - 1$ and $\sg_b = 2 b$. Again we study two cases as $8 n - 7$ can be a perfect square or not.

\begin{enumerate}
\item If $8 n - 7$ is a perfect square, then there exists such integer $m$ that $8 n - 7 = m^2$. Thus $m$ is odd and $m \geq 5$ as $n \geq 6$ (and thus $8 n - 7 \geq 35 > 5^2$). Hence  $$a = \left \lceil \frac{1 + m}{2} \right \rceil = \frac{1 + m}{2} \; \text{ and } \; \sg_a = m.$$

From $8 n - 7 = m^2$ and $m \geq 4$ it follows that $m^2 < 8 n + 1 \leq (m+1)^2$.
Thus
$$b = \left \lceil \frac{- 1 + (m + 1)}{2} \right \rceil = \left \lceil \frac{m}{2} \right \rceil = \frac{m+1}{2} \; \text{ and } \; \sg_b = m + 1.$$

Hence $\sg_a < \sg_b$ and the optimal value is reached for $s_2 = s_1 - 1$. Additionally we notice that $a = b$ so $\sg(K_{n, n}) = 2 b - 1$.

\item If $8 n - 7$ is not a perfect square, then there exists such integer $m$ that $$m^2 < 8n - 7 < (m + 1)^2.$$ Thus $$a = \left \lceil \frac{1 + \sqrt{8 n - 7}}{2} \right \rceil = \left \lceil \frac{1 + (m+1)}{2} \right \rceil = \left \lceil \frac{m + 2}{2} \right \rceil$$
and
$$\sg_a = \begin{cases}
m + 1; & m \text{ even},\\
m + 2; & m \text{ odd}.
\end{cases}$$

As $m \geq 3$, it also holds that $8 n + 1 < (m + 1)^2 + 8 < (m + 2)^2$. Hence
$$b = \left \lceil \frac{- 1 + \sqrt{8 n + 1}}{2} \right \rceil \leq \left \lceil \frac{- 1 + (m + 2)}{2} \right \rceil = \left \lceil \frac{m + 1}{2} \right \rceil.$$

If $m$ is odd, then $$\sg_b = 2 \cdot \frac{m+1}{2} = m+1 < m+2 = \sg_a,$$
thus the optimal value is reached for $s_1 = s_2$. 

If $m$ is even, then $m+1$ is odd. By Lemma~\ref{lihi}, there exists an integer $k$ such that $8 k + 1 = (m+1)^2$. Thus $(m+1)^2 \equiv 1 \pmod 8$ and also $8 n - 7 \equiv 1 \pmod 8$. Hence  $8 n - 7 \notin \{(m+1)^2 - 7, \ldots, (m+1)^2- 1, (m+1)^2 + 1, \ldots, (m+1)^2 + 7\}$. This means that the case $8 n + 1 > (m+1)^2$ is not possible. Thus we can improve the bound for $b$:
$$b = \left \lceil \frac{- 1 + \sqrt{8 n + 1}}{2} \right \rceil \leq \left \lceil \frac{- 1 + (m + 1)}{2} \right \rceil = \left \lceil \frac{m}{2} \right \rceil =  \frac{m}{2}.$$
Hence $$\sg_b = m < m + 1 = \sg_a.$$ The optimal value is reached for $s_1 = s_2$. It follows from the above that if $8 n - 7$ is not a perfect square, we have $\sg(K_{n,n}) = 2 b$. \qed
\end{enumerate}

To conclude the section we determine $\sg(K_{n_1, n_2})$ for the case when $n_1 \gg n_2$. More precisely:

\begin{proposition}
If $n_1, n_2 \in \N, n_1, n_2 \geq 2$, and $\displaystyle \binom{n_1 - \binom{n_2}{2}}{2} \geq n_1$, then
$$\sg(K_{n_1, n_2}) = 
\begin{cases}
n_1 ; & n_2 = 2,\\
n_1 + n_2 - \binom{n_2}{2}; & n_2 \geq 3.
\end{cases}$$
\end{proposition}

\proof
Let $(X, Y)$ be a bipartition of the graph $K_{n_1, n_2}$ and $S = S_1 \cup S_2$, $S_1 \subseteq X, S_2 \subseteq Y$, its strong geodetic set. Let $|S_i| = s_i$ for $i \in [2]$. Hence $\sg(K_{n_1, n_2}) = s_1 + s_2$.

Suppose $n_2 = 2$, then $s_2 \in \{0, 1, 2\}$. For each of these three possibilities we consider the possible values of $s_1$ (cf.~Table~\ref{tabela}), such that all the vertices of the graph are covered.

\begin{table}[!h]
\begin{center}
\begin{tabular}{c c | c}
$s_2$ & $s_1$ & $\sg$ \\ \hline
$0$ & $n_1$ & $n_1$ \\
$1$ & $n_1$ & $n_1 + 1$ \\
$2$ & $n_1 - 1$ & $n_1 + 1$\\
\end{tabular}
\end{center}
\caption{Possible values of $s_1$ and $s_2$ in the case when $n_2 = 2$.}
\label{tabela}
\end{table}


Clearly, the smallest value is reached for $S = X$. Hence we have $\sg(K_{n_1, 2}) = n_1$.

Suppose $n_2 \geq 3$. Thus $\binom{n_2}{2} \geq n_2$ (which does not hold for the above case). In addition, the sequence $\{\binom{k}{2} - k\}_{k \in \N, k \geq 3}$ is increasing.

If the strong geodetic set contains $s_2$ vertices from $Y$, then these vertices can cover $\binom{s_2}{2}$ vertices in $X$. From $\displaystyle \binom{n_1 - \binom{n_2}{2}}{2} \geq n_1$ it follows that the remaining $s_1 = n_1 - \binom{s_2}{2}$ vertices from $X$ can cover all still uncovered vertices in $Y$. Hence the optimization problem~\eqref{optimizationProblem} can be simplified to:
\begin{align*}
\min \quad & n_1  + s_2 - \binom{s_2}{2} \\
\text{subject to: } & 0 \leq s_2 \leq n_2\\
& s_2 \in \Z.
\end{align*}
As the sequence $\{\binom{k}{2} - k\}$ is increasing, the sequence $\left \{n_1 - \left(\binom{k}{2} - k \right) \right \}$ is decreasing and its minimum is reached at the largest possible value of $k$. The solution of this optimization problem is $s_2 = n_2$. Hence $\sg(K_{n_1, n_2}) = s_1 + s_2 = n_1 + n_2 - \binom{n_2}{2}$.
\qed

\begin{remark}
The above proposition holds (with the same proof) in a more general case when for all values of $s_2$ such that $0 \leq s_2 \leq n_2$, it holds $ \displaystyle\binom{n_1 - \binom{s_2}{2}}{2} \geq n_1 - s_2$.
\end{remark}

\section{The strong geodetic number versus diameter}
\label{sec:diameter}

Chartrand et al.\ \cite{obstoj} studied the geodetic number of a graph with specified diameter. In this section we derive similar results for the strong geodetic number. 

Taking into account the diameter of a graph we can improve Proposition~\ref{prop:basic} as follows.

\begin{proposition}
\label{prop:diameter}
If $G$ is a graph with $\n(G) \geq 2$, then $$\sg(G) \leq \n(G) - \diam(G) + 1.$$
\end{proposition}

\proof
Let $n = \n(G)$, $d = \diam(G)$ and $u, v \in V(G)$ such that $\d(u, v) = d$. Let $u = v_0, v_1, \ldots, v_{d-1}, v_d = v$ be a path $P$ of length $d$ between $u$ and $v$. Define a set $S = V(G) - \{v_1, \ldots, v_{d-1}\}$. Then it is clear that $S$ is a strong geodetic set and therefore, $\sg(G) \leq n - d + 1$.
\qed

The example showing that the inequality in Proposition~\ref{prop:diameter} is best possible is the same graph that attains the equality in $\g(G) \leq \n(G) - \diam(G) + 1$ in~\cite{obstoj}.

Using a graph from~\cite{obstoj} with raduis $r$, diameter $d$ and geodetic number $k$ we derive the following result (but only for $k \geq 3$).

\begin{proposition}
\label{prop:radius}
If $r, d, k \in \mathbb{N}$, $k \geq 3$ and $r \leq d \leq 2 r$, then there exists a connected graph $G$ such that $\rad(G) = r$, $\diam(G) = d$ and $\sg(G) = k$.
\end{proposition}

\begin{remark}
Proposition~\ref{prop:radius} does not hold for $k = 2$. From $\sg(G) = 2$ it follows that $G$ is a path, thus $\rad(G) = \left \lfloor \frac{\diam(G)}{2} \right \rfloor$ and the radius may not be chosen arbitrarily.
\end{remark}

\medskip

In~\cite{Fitzpatrick-1999-CN} we find the following property of the isometric path number, which can be generalized to the strong geodetic number. For a graph $G$ it holds
$$\ip(G) \geq \frac{\n(G)}{\diam(G) + 1}.$$
It is stated in~\cite{Fitzpatrick-1999-CN} that the equality is reached for paths, cycles, complete graphs and hypercubes on $2^{2^k - 1}$ vertices, which is not entirely correct. The equality holds for paths as $\ip(P_n) = 1 = \frac{n}{(n-1) + 1}$ and for complete graphs of even order as $\ip(K_{2k}) = k = \frac{2k}{2}$. But for cycles we have $\ip(C_n) = 2 > \frac{n}{\lfloor \frac{n}{2} \rfloor + 1}$, thus the equality is never reached. Similarly it does not hold for complete graphs of odd order as $\ip(K_{2k+1}) = k+1 >\frac{2k+1}{2}$. The equality in $\ip(G) \geq \left \lceil \frac{\n(G)}{\diam(G) + 1} \right \rceil$ is however attained also for cycles and odd complete graphs.

We now consider the analogous relation for strong geodetic number.

\begin{proposition}
\label{prop:first}
If $G$ is a graph with $\n(G) \geq 2$, then 
$$\sg(G) \geq \left \lceil \frac{1 + \sqrt{1 + \frac{8 \n(G)}{\diam(G) + 1}}}{2} \right \rceil.$$
\end{proposition}

\proof
Let $n = \n(G)$ and $d = \diam(G)$. As the diameter of the graph $G$ equals $d$, each geodesic has length at most $d$ and thus covers at most $d+1$ vertices. As the graph is covered with $\binom{\sg(G)}{2}$ geodesics, it follows that
$$n \leq \binom{\sg(G)}{2} (d+1),$$ 
which in turn implies 
$$\sg(G)^2 - \sg(G) - \frac{2n}{d+1} \geq 0.$$
One of the corresponding zeros is negative and since $\sg(G)$ is an integer, we conclude that $$\sg(G) \geq \left \lceil \frac{1 + \sqrt{1 + \frac{8 n}{d + 1}}}{2} \right \rceil.$$
\qed

Equality is attained for $K_1$ and $P_n$ ($n \geq 2$) as $d = n-1$ and $\sg(P_n) = 2$. Furthermore, equality can only be reached if every vertex (including endpoints of geodesics) is covered exactly once. Thus $\sg(G) \leq 2$, which implies $\sg(G) = 2$, hence $G$ must be isomorphic to a path.

Observe that in Proposition~\ref{prop:first} every vertex from the strong geodetic set is covered $\sg(G)$-times. Thus we can develop a better bound for $\sg(G)$. But before that we consider a special case when the diameter of a graph equals $1$. Only such graphs are complete graphs and for them we know that $\sg(G) = n$. 

The main result of this section reads as follows.

\begin{theorem}
\label{thm:second}
If $G$ is a graph with $\n(G) = n$ and $\diam(G) = d \geq 2$, then 
$$\sg(G) \geq \left \lceil \frac{d - 3 + \sqrt{(d-3)^2 + 8 n (d-1)}}{2 (d-1)} \right \rceil.$$
\end{theorem}

The proof uses the observation $\n(G) \leq \sg(G) + \binom{\sg(G)}{2} (\diam(G) - 1)$ and is similar to the proof of Proposition~\ref{prop:first}. The equality is attained for example for paths. We now consider some other equality cases.

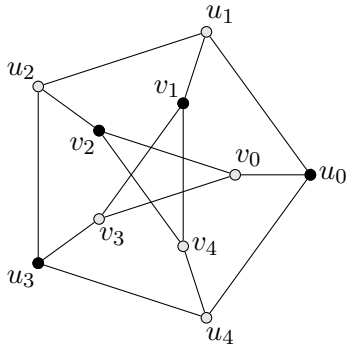
\begin{figure}[!ht]
\begin{center}
    \begin{tikzpicture}[scale=0.5]
    \pgfmathtruncatemacro{\N}{5}
    \pgfmathtruncatemacro{\Nr}{4}
    \pgfmathtruncatemacro{\r}{2}
    \pgfmathtruncatemacro{\R}{4}
    \begin{scope}
        \foreach \x in {0,...,\Nr}
            \node[label={\x*360/\N:$u_{\x}$}] (\x) at (\x*360/\N:\R cm) {};
        \foreach \x [remember=\x as \lastx (initially 0)] in {0,...,\Nr,0}
            \path (\x) edge (\lastx);
    \end{scope}
    
    \foreach \x  in {0, 3}
        \node[fill=black] () at (\x*360/\N:\R cm) {};
        
    \begin{scope}
        \foreach \x in {0,...,\Nr}
            \node[label={(\x+1)*360/\N:$v_{\x}$}] (\x+5) at (\x*360/\N:\r cm) {};
    \end{scope}
    
    \foreach \x  in {6, 7}
        \node[fill=black] () at (\x*360/\N:\r cm) {};
        
    \draw (0+5) edge (2+5); 
    \path (2+5) edge (4+5); 
    \path (4+5) edge (1+5); 
    \path (1+5) edge (3+5); 
    \path (3+5) edge (0+5);
    
    \foreach \x  in {0,...,\Nr}
        \path (\x) edge (\x+5);
    
    \end{tikzpicture}
    \caption{A strong geodetic set of the Petersen graph.}
    \label{petersen}
\end{center}
\end{figure}

Let $P$ be the Petersen graph with vertex labeling as shown in Fig.~\ref{petersen}. It follows from Theorem~\ref{thm:second} that $\sg(P) \geq 4$. The set $S = \{ u_0, u_3, v_1, v_2 \}$ is a strong geodetic set, so the equality is attained.

\begin{figure}[!ht]
\begin{center}
    \begin{tikzpicture}[scale=0.6]
    \pgfmathtruncatemacro{\R}{4}
    \begin{scope}
        \foreach \x in {1,...,3}
            \node[fill=black] (\x) at (\x*360/3:\R cm) {};
        \foreach \x [remember=\x as \lastx (initially 1)] in {1,...,3,1}
            \path (\x) edge (\lastx);
    \end{scope}
    
    \node[fill=black] (0) at (0, 0) {};

    \draw (0) edge (1); 
    \path (0) edge (2); 
    \path (0) edge (3); 
    
    \node (a) at (1.5*360/3:2 cm) {};
    \node (b) at (2.5*360/3:2 cm) {};
    \node (c) at (0.5*360/3:2 cm) {};
    \node (a') at (1*360/3:2 cm) {};
    \node (b') at (2*360/3:2 cm) {};
    \node (c') at (3*360/3:2 cm) {};
    
    \draw (a) edge (b);
    \draw (b) edge (c);
    \draw (c) edge (a);
    \draw (a') edge (b');
    \draw (b') edge (c');
    \draw (c') edge (a');
    \draw (a) edge (a');
    \draw (a') edge(c);
    \draw (a) edge (b');
    \draw (b) edge (b');
    \draw (b) edge (c');
    \draw (c) edge (c');

    \draw [bend left=20,-] (b) to (a');
    \draw [bend left=20,-] (c) to (b');
    \draw [bend left=20,-] (a) to (c');
    
    \end{tikzpicture}
    \caption{The graph $G_4$ with strong geodetic number $4$ and diameter $2$ which attains the equality in Theorem~\ref{thm:second}.}
    \label{subdivizija}
\end{center}
\end{figure}
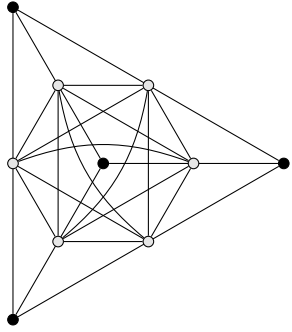

We now consider more general graphs of diameter $2$ which attain the equality in Theorem~\ref{thm:second}. Let $G$ be a graph of diameter $2$. Let $S$ be its strong geodetic set of size $\sg(G)$. Thus $\n(G) = \sg(G) + \binom{\sg(G)}{2}$ and all geodesics must be of length $2$. Hence vertices of graph $G$ can be partitioned in two parts: vertices lying in $S$ and vertices which lie in exactly one geodesic between vertices in $S$. The set $S$ is an independent set, otherwise the geodesic between connected vertices would not be of length $2$. Consider the following construction of a graph $G_k$ for an integer $k$. Take a subdivision of $K_{k}$ with vertices $V(K_k) \cup E(K_k)$ and add edges $e \sim f$ for each pair $e, f \in E(G)$  (cf.~Fig.~\ref{subdivizija}). Clearly, the diameter of $G_k$ is $2$. The set $V(K_k)$ with subdivided edges as fixed geodesics is a strong geodetic set, thus $\sg(G_k) \leq k$. But it follows from $\n(G_k) = k + \binom{k}{2}$ and Theorem~\ref{thm:second} that $k \leq \sg(G_k)$. Hence, $\sg(G_k) = k$ and the equality in Theorem~\ref{thm:second} is attained for $G_k$.  

The above observations can be generalized to an arbitrary diameter. Let $G$ be a graph with diameter $d$, strong geodetic number $k$ and $\n(G) = k + (d-1) \binom{k}{2}$. If $S$ is a strong geodetic set of size $k$, then it is independent and between each pair of vertices from $S$, there exist completely disjoint paths of length $d$ (thus each path contains $d-1$ vertices). These paths are geodesics, so there are no other edges between vertices on the same path. However, edges between vertices on different paths must be added in such a way that $G$ has diameter $d$. 

\begin{figure}[!ht]
\begin{center}
    \begin{tikzpicture}[scale=0.7]
    \pgfmathtruncatemacro{\R}{4}
    \begin{scope}
        \foreach \x in {1,...,3}
            \node[fill=black] (\x) at (\x*360/3:\R cm) {};
        \foreach \x [remember=\x as \lastx (initially 1)] in {1,...,3,1}
            \path (\x) edge (\lastx);
    
    \node[fill=black] (0) at (0, 0) {};

    \draw (0) edge (1); 
    \path (0) edge (2); 
    \path (0) edge (3); 
    
    \node (a) at (1.5*360/3:2 cm) {};
    \node (b) at (2.5*360/3:2 cm) {};
    \node (c) at (0.5*360/3:2 cm) {};
    \node (a') at (1*360/3:2 cm) {};
    \node (b') at (2*360/3:2 cm) {};
    \node (c') at (3*360/3:2 cm) {};
    \node (a1) at (1.1*360/3:3 cm) {};
    \node (a2) at (1.9*360/3:3 cm) {};
    \node (b1) at (2.1*360/3:3 cm) {};
    \node (b2) at (2.9*360/3:3 cm) {};
    \node (c1) at (3.1*360/3:3 cm) {};
    \node (c2) at (3.9*360/3:3 cm) {};
    \node (a1') at (1*360/3:3 cm) {};
    \node (a2') at (1*360/3:1 cm) {};
    \node (b1') at (2*360/3:3 cm) {};
    \node (b2') at (2*360/3:1 cm) {};
    \node (c1') at (3*360/3:3 cm) {};
    \node (c2') at (3*360/3:1 cm) {};
    
    \draw [bend right=10,-] (a) to (b);
    \draw [bend right=10,-] (b) to (c);
    \draw [bend right=10,-] (c) to (a);
    \draw (a') edge (b');
    \draw (b') edge (c');
    \draw (c') edge (a');
    \draw (a) edge (a');
    \draw (a') edge(c);
    \draw (a) edge (b');
    \draw (b) edge (b');
    \draw (b) edge (c');
    \draw (c) edge (c');

    \draw [bend left=20,-] (b) to (a');
    \draw [bend left=20,-] (c) to (b');
    \draw [bend left=20,-] (a) to (c');
    \end{scope}
    
    \begin{scope}[xshift=9cm]
        \foreach \x in {1,...,3}
            \node[fill=black] (\x) at (\x*360/3:\R cm) {};
        \foreach \x [remember=\x as \lastx (initially 1)] in {1,...,3,1}
            \path (\x) edge (\lastx);
    
    \node[fill=black] (0) at (0, 0) {};

    \draw (0) edge (1); 
    \path (0) edge (2); 
    \path (0) edge (3); 
    
    \node (a) at (1.75*360/3:2.3 cm) {};
    \node (b) at (2.75*360/3:2.3 cm) {};
    \node (c) at (0.25*360/3:2.3 cm) {};
    \node (a1) at (1.25*360/3:2.3 cm) {};
    \node (b1) at (2.25*360/3:2.3 cm) {};
    \node (c1) at (0.75*360/3:2.3 cm) {};
    \node (a') at (1*360/3:1 cm) {};
    \node (b') at (2*360/3:1 cm) {};
    \node (c') at (3*360/3:1 cm) {};
    \node (a1') at (1*360/3:2.3 cm) {};
    \node (b1') at (2*360/3:2.3 cm) {};
    \node (c1') at (3*360/3:2.3 cm) {};
    
    \draw (a) edge (b);
    \draw (b) edge (c1);
    \draw (c1) edge (a);
    \draw (a') edge (b');
    \draw (b') edge (c');
    \draw (c') edge (a');
    \draw (a) edge (a');
    \draw (a') edge(c1);
    \draw (a) edge (b');
    \draw (b) edge (b');
    \draw (b) edge (c');
    \draw (c1) edge (c');

    \draw (b) edge (a');
    \draw (c1) edge (b');
    \draw (a) edge (c');
    \end{scope}
    
    \end{tikzpicture}
    \caption{Graphs $G_{4, 4}$ and $G_{4, 3}$ both attain equality in Theorem~\ref{thm:second}.}
    \label{subdivizijaPremer}
\end{center}
\end{figure}
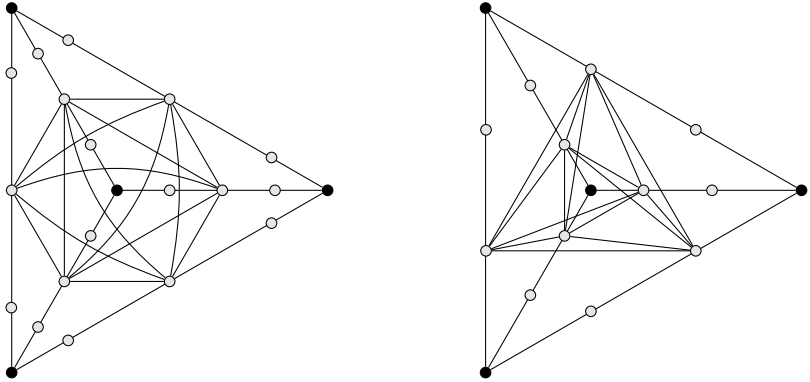

Consider the following graph $G_{k, d}$. Take a $(d-1)$-times subdivision of a graph $K_k$ with vertices $V(K_k) \cup \{ e^i \; ; \; e \in E(K_k), i \in [d-1] \}$, where $e^i$ all lie on a subdivided edge $e$ and $e^i \sim e^{i+1}$ for all $i \in [d-2]$. Define a set $\V \subseteq V(G_{k, d})$ such that for each $u \in V(K_k)$ there exists a vertex $v \in \V$ such that $\d(u, v) = \left \lfloor \frac{d}{2} \right \rfloor$. Add edges $x \sim y$ to the graph $G_{k, d}$ for all $x, y \in \V$, $x \neq y$ (cf.~Fig.~\ref{subdivizijaPremer}). Observe that when $d$ is even, the set $\V$ is unique and equals $\{ e^{\frac{d}{2}} \; ; \; e \in E(K_k) \}$. 

Clearly, $\diam(G_{k, d}) \geq d$, as $\d(u, v) = d$ for any $u, v \in V(K_k)$. First suppose $d$ is even. For each $e^i \in V(G_{k, d}) - V(K_k)$ there exists $v \in \V$ such that $\d(e^i, v) \leq \frac{d}{2} - 1$. Thus $\d(e^i, u) \leq \frac{d}{2} - 1 + 1 + \frac{d}{2} = d$ for any $u \in V(G_{k, d})$. So in this case we have $\diam(G_{k, d}) = d$. 

Now suppose $d$ is odd. For each $e^i \in V(G_{k, d}) - V(K_k)$ there exists $v \in \V$ such that $\d(e^i, v) \leq \frac{d-1}{2}$. Thus for any $e^i, f^j \in V(G_{k, d}) - V(K_k)$ it holds $\d(e^i, f^j) \leq \frac{d-1}{2} + 1 + \frac{d-1}{2} = d$. Due to the property of the set $\V$, between vertices $e^i \in V(G_{k, d}) - V(K_k)$ and $u \in V(K_k)$ there exists a path of length at most $\frac{d-1}{2} + 1 + \left \lfloor \frac{d}{2} \right \rfloor = d$. Hence, $\diam(G_{k, d}) = d$ in both cases.

As the set $V(K_k)$ is a strong geodetic set and $\n(G_{k, d}) = k + (d-1) \binom{k}{2}$, it holds $\sg(G_{k, d}) = k$. Therefore, the inequality in Theorem~\ref{thm:second} is attained for $G_{k, d}$. 

\section{Further research}
\label{sec:Fur-Research}

In this paper we have studied the strong geodetic problem for complete bipartite graphs. The first natural problem is to determine the strong geodetic number for all $K_{n_1, n_2}$ and also to consider complete multipartite graphs. Discovering some more examples of graphs which attain the equality in Theorem~\ref{thm:second} would also be of great interest, and ideally, to characterize such graphs.

%

\end{document}